\documentclass[10pt,twocolumn,twoside]{IEEEtran} 

\IEEEoverridecommandlockouts                              

\usepackage{amsmath,amsfonts,amssymb,bm}
\usepackage{multicol,multirow,diagbox,graphicx,subfigure,lineno,cases,cite}
\usepackage{epstopdf,pgf,tikz,color,xcolor}
\usepackage{graphicx,subfigure,pgf,tikz,pst-node}
\usetikzlibrary{arrows,shapes,automata,positioning,fit}
\allowdisplaybreaks 

\newcommand{\e}{\varepsilon}
\newcommand{\w}{\omega}
\newcommand{\N}{\mathcal{N}}
\newcommand{\G}{\mathcal{G}}

\renewcommand{\bar }{\overline}
\newcommand{\R}{\mathbb{R}}
\newtheorem{ass}{Assumption}
\newtheorem{lem}{Lemma}
\newtheorem{thm}{Theorem}
\newtheorem{rem}{Remark}
\newtheorem{prob}{Problem}

\renewcommand{\bar }{\overline}
\newcommand{\pb}{\begin{IEEEproof} }
\newcommand{\pe}{\end{IEEEproof}}

\begin{document}
	\title{On Measurement Disturbances in Distributed Least Squares Solvers for Linear Equations }
	\author{Yutao Tang, Yichen Zhang, Ruonan Li, and Xinghu Wang 
		\thanks{This work has been submitted to the IEEE for possible publication. Copyright may be transferred without notice, after which this version may no longer be accessible.}
		\thanks{This work was supported by National Natural Science Foundation of China under Grant 61973043.}
		\thanks{Y. Tang, Y. Zhang, and R. Li are with the School of Artificial Intelligence, Beijing University of Posts and Telecommunications, Beijing 100876, China (e-mails: yttang@bupt.edu.cn, zhangyc930@163.com, nanruoliy@163.com). X. Wang is with Department of Automation, University of Science and Technology of China, Hefei, 230027, China (e-mail: xinghuw@ustc.edu.cn). 
	}
	}

\maketitle

\begin{abstract}
	This paper aims at distributed algorithms for solving a system of linear algebraic equations. Different from most existing formulations for this problem, we assume that the local data at each node is not accurately measured but subject to some disturbances. To be specific,  the local measurement consists of two parts: a nominal value and a multiple sinusoidal disturbance. By introducing an identifier-enhanced observer to estimate the disturbance, we present a novel distributed least squares solver for the linear equations using noisy measurements. The proposed solver is proven to be able to recover the least squares solution to the linear equations associated with the nominal values irrespective of any multi-sinusoidal disturbance even with unknown frequencies. We also show the robustness of the distributed solvers under standard conditions against unstructured perturbations. The effectiveness of our design is verified by a numerical example.
\end{abstract}


\section{Introduction}\label{sec:intro}

In this paper, we are concerned with the following linear algebraic equations:
\begin{align}\label{eq:main}
	z=Hy
\end{align}
where $y \in \mathbb{R}^m$,  $z(t)\in \mathbb{R}^M$,  and $H \in \mathbb{R}^{M\times m}$.  This system of linear equations might be the most well-known equation in algebra and has been extensively used to describe or approximate the relations between two different variables  arising in natural science, engineering, and social sciences \cite{boyd2018introduction, strang2019linear}.  
Over the past few years, this problem involving large-scale networks has gained renewed interests due to the fast development of sensor network and big data technologies. In this paper, we aim at distributed continuous-time solvers to recover the least squares solution to \eqref{eq:main} using noisy measurements. 


{\em Literature review}:  As a fundamental problem in scientific computing, this problem has been well-studied in literature and efficient centralized algorithms can be found in classical textbooks, e.g., \cite{golub2013matrix}.  To solve this problem in a distributed manner by a group of interconnected compute nodes (or agents), each node typically knows a part of the matrix $H$ and thus partial measurement $z$, and cooperatively maintains an estimate about the global solution. These nodes share its information with other nodes and iteratively update its own estimate. In an early attempt \cite{pasqualetti2012distributed}, each node should share its estimate and certain information of $H$ (e.g., the known rows or its corresponding kernel space). In the seminal work \cite{mou2015distributed}, a constrained consensus-based algorithm was constructed to solve the problem where each node only needs to share its own estimate of the solution under local feasibility initialization requirement. Similar algorithms were proposed in \cite{shi2017network, liu2017exponential}. To remove the local initialization requirement,  distributed optimization-based algorithms have attracted much attention and been modified to solve the linear equations in the literature, e.g., \cite{liu2017asynchronous, wang2019initialization, wang2019distributed, liu2020network} and references therein. Recently, such optimization-based designs have been extended to solve more general matrix equations in \cite{zeng2018distributed, deng2022network}. Note that a large portion of aforementioned results highly rely the perfect measurement of $z$ to solve the problem. 

In practical scenarios,  the data are typically collected from sensor networks and inevitably affected by some measurement disturbance/noise. On the other hand, even with accurate measurements, the algorithms might suffer from computation errors during its implementation. Thus, it is of practical interest to develop effective algorithms that can deal with such disturbances arising from either inaccurate measurements or computation errors. 

A popular treatment in the literature is to model the disturbance as a zero-mean random noise and formulate the problem as a linear filtering problem. 
Another treatment is to view the disturbance as some unstructured but bounded perturbation and pursue the worst case solution in a robust sense. Both treatments have been conducted to solve linear equations with measurement uncertainties in centralized settings \cite{sayed2011adaptive, nagpal1991filtering} with a few very recent distributed extensions, e.g., \cite{yang2020distributed, yi2022distributed, lee2020distributed}. In these existing results, the useful part of the full measurement is assumed to be a constant while the rest part is roughly viewed as zero-mean stochastic or unstructured deterministic disturbances. However, in practice, the measurement disturbance can be time-varying and heavily deteriorate the algorithm performance. As such, how to actively reject such nontrivial measurement disturbances rather than passively accepting the adverse effect caused by them should be carefully considered in the design of distributed linear equation solvers using noisy measurements.  

Meanwhile, there are some interesting attempts on the disturbance rejection issue in distributed optimization algorithms. For example, the authors combined the distributed gradient dynamics with internal models and solve distributed optimization problem in spite of sinusoidal actuated disturbance in \cite{wangxinghu2017}. In \cite{tang2018smc},  an observer-based design was presented to actively estimate such external disturbances and ensure the distributed optimization goal with disturbance rejection.  Similar designs can be found in \cite{wangzhaojian2019} to handle unknown time-varying power imbalance in distributed load-side frequency control. However, the considered disturbances in these results are assumed in the actuator channel rather than measurement channel.  Thus, the aforementioned distributed optimization algorithms cannot be directly utilized to solve the linear equations and reject possible measurement disturbances. 





{\em Statement of contributions}: 
We propose a novel distributed least squares solver for the given linear equations using only noisy measurements. The main design is a novel combination of distributed primal-dual dynamics and identifier-based disturbance compensator. We prove that the presented solver can effectively recover the expected least squares solution under weight-balanced digraphs irrespective of any multi-sinusoidal measurement disturbance. 
Moreover, the convergence rate is shown to be exponentially fast from any initial point. We also show the robustness of our designs against unstructured perturbations. Numerical simulations illustrate our results.


{\em Organization}: The rest of this paper is organized as follows. We first state the problem in Section \ref{sec:form}. Then, we constructively present a novel distributed solver under weight-balanced digraphs and then enhance it with a frequency identifier to solve our formulated problem in Section \ref{sec:alg}. In Section \ref{sec:sol}, we prove the effectiveness of proposed algorithms to solve the problem with disturbance rejection and also confirm their robustness with respect to unstructured perturbations. A simulation example is given in Section \ref{sec:simu} to validate our design.  Section \ref{sec:con} closes the paper with some brief remarks.


\section{Problem statement}\label{sec:form}

Consider the system of linear equations \eqref{prob}. We are interested in distributed solvers to recover its least squares solution formulated as a minimizer to the following optimization problem:
\begin{align}\label{prolem:opt}
	{y}^*\in \arg\min  f({y})\triangleq || z-{H} {y}||^2 
\end{align}
To be specific, we consider $N$ interconnected compute nodes. Each node $i$ knows $N_i$ mutually different rows of ${H}$ and the corresponding portion of the full measurement: 
\begin{align}\label{eq:main-noise}
	\tilde z_i={H}_i{y}+\e_i(t)
\end{align}
where $z_i=H_i y\in \R^{N_i}$ represents the perfect measurement (or nominal value) at node $i$ and $\e_i(t) \in \R^{N_i}$ is the local measurement disturbance with $H=[H_1^\top,\,\dots,\,H_N^\top ]^\top$.  Let $\e=[\e_1^\top,\,\dots,\,\e_{N}^\top]^\top$ for short.  Here we suppose $\e(t)$ may heavily distort the  nominal value and focus on effective solvers for the nominal equations \eqref{eq:main} with disturbance rejection.

Note that the optimization problem \eqref{prolem:opt} is convex and always admits at least one global solution. Without loss of generality, we suppose the matrix ${H}$ has a full column rank, i.e., ${\rm rank}(H)=m$. Then the least squares solution to \eqref{eq:main} 
is unique and can be expressed as  $y^*=(H^\top H)^{-1}H^\top z$. 

Suppose node $i$ maintains an estimate $x_i\in \R^m$ of the expected least squares solution $y^*$ and can share its information with its immediate neighbors through an underlying communication topology described by a digraph $\G=\{\mathcal{N},\,\mathcal{E},\,A\}$ with node set $\N=\{1, \,\dots,\,N\}$, edge set $\cal E\subset \N\times \N$, and the adjacency matrix $A=[a_{ij}]_{N\times N}$. An edge $(i,\,j)\in \cal E$ means node $j$ can get access to the information of node $i$. 




It is ready for us to formulate our problem as follows: 

\begin{prob}\label{prob}
	For a given digraph $\G$, matrix ${H}_i$, and $\tilde z_i(t)$, find an effective distributed update law for each node relying only local communication and computation 
	such that, for any $x_i(0)\in \R^m$ and $\e_i(t)\in \R^{m}$,  the trajectory $x_i(t)$ is uniformly bounded and satisfies that $\lim_{t\to\infty} x_i(t)={y}^*$.
\end{prob}

When $\e(t)\equiv {\bf 0}$,  this formulation exactly coincides with the linear least squares problem for the equations \eqref{eq:main} discussed in the literature \cite{mou2015distributed, shi2017network, liu2017exponential, yang2020distributed,zeng2018distributed}. With this measurement disturbance, we attempt to develop an effective mechanism to asymptotically solve the corresponding nominal liner equations \eqref{eq:main} with disturbance rejection. 
As a result, our problem is more technically challenging than the aforementioned results where the measurement disturbance is not considered or passively treated as zero-mean stochastic noises or unstructured bounded perturbations in existing works.

To ensure the tractability of our problem,  we impose the following assumptions on the communication graph and the generic structure of measurement disturbances. 

\begin{ass}\label{ass:graph}
	Digraph $\G$ is weight-balanced and strongly connected. 
\end{ass}

\begin{ass}\label{ass:disturbance}
	For any $i\in \N$,  each component of $\e_i(t)$ contains only $k_i$ sinusoidal signals with distinct frequencies $\omega_{i1}$, $\dots$, $\omega_{i k_i}$. 
\end{ass}

Assumption \ref{ass:graph} is about the connectivity of communication graphs for these nodes. 
Assumption \ref{ass:disturbance} specifies the class of measurement disturbances we are interested in. In practice, we can truncate the dominated finite harmonics of any bounded periodic signal in its Fourier series expansion and have a fair approximation fulfilling this assumption. In this sense, this assumption is not restrictive as it appears. We will discuss the robustness of our algorithms with respect to unstructured disturbances later.  Both assumptions have been widely used in the multi-agent coordination or disturbance rejection literature \cite{shi2017network, liu2017exponential, zeng2018distributed, marino2003output, sastry2011adaptive, tang2018iet, wangxinghu2017, wangzhaojian2019}.

From the expression of $y^*$, the expected least squares solution corresponding to the nominal equations \eqref{eq:main} will be recovered when we have a good enough estimation of the nonvanishing measurement disturbance  $\e(t)$. Thus, the disturbance rejection property is essentially equivalent to solve an estimation or signal reproduction problem. For the ease of presentation, we suppose $N_i=1$ throughout this paper and present a novel adaptive observer-based design to solve the formulated problem via noisy measurements.

\section{Algorithm Design}\label{sec:alg}

In this section, we present our main design to solve the linear equations with disturbance rejection.  We will first develop a distributed solver for the case with accurate measurements, and then present an identifier-enhanced distributed algorithm to solve the formulated problem with disturbance rejection. 

\subsection{Distributed solver with accurate measurements}

We suppose the accurate measurements are available without any  disturbance, i.e, $\tilde z_i=z_i$. The problem becomes the conventional form  discussed in the literature \cite{mou2015distributed, wang2019solving}. 

Under Assumption \ref{ass:graph}, we let $\mbox{Sym}(L)=\frac{L+L^\top}{2}$ with $L$ the Laplacian of $\G$.  It can be verified that $\mbox{Sym}(L)$ is positive semidefinite with $0$ as its simple eigenvalue. Then we can order the eigenvalues of $\mbox{Sym}(L)$ as $0=\lambda_1<\lambda_2\leq \dots\leq \lambda_N$. Choose a matrix ${R}\in \R^{N\times (N-1)}$ such that ${R}^\top {\bf 1}_N={\bf 0}$, ${R}^\top {R}=I_{N-1}$, and ${R}{R}^\top+\frac{{\bf 1}_N {\bf 1}_N^\top }{N}=I_m$. One can further verify that $ \lambda_2 I_{N-1}\leq R^\top \mbox{Sym}(L)R\leq \lambda_N I_{N-1}$ holds under Assumption \ref{ass:graph}.

Motivated by the distributed optimization designs in \cite{wang2011control, gharesifard2014distributed, tang2019optimal, liu2019arrow}, we present the following  algorithm to solve \eqref{eq:main}:
\begin{align}\label{alg:full-exact}
	\dot{x}_i&=- \kappa_1 {H}_i^\top [{H}_ix_i-z_i]-(\kappa_1+\kappa_2) x_{oi}-v_{oi}\\
	\dot{v}_i&=\kappa_1\kappa_2 x_{oi}
\end{align}
with $x_{oi}=\sum_{j=1}^N a_{ij}(x_i-x_j)$ and $v_{oi}=\sum_{j=1}^N a_{ij}(v_i-v_j)$. Here $v_i\in \R^m$ is an auxiliary variable and $\kappa_1,\,\kappa_2>0$ are parameters to be specified later.  It is noted that \eqref{alg:full-exact} is in the form of primal-dual dynamics to solve a distributed optimization problem with local cost function $f_i(y)=\frac{1}{2}||{H}_iy-  z_i||^2$. Nevertheless, since $f_i$ is only convex but not strongly convex,  the proofs in \cite{kia2015distributed, tang2021optimal} cannot be directly implemented to show the exponential convergence of \eqref{alg:full-exact}.

Denote by ${\lambda}_{h}$ and $ {\lambda}_H $ the minimal and maximal eigenvalue of matrix $H^\top H$ and let $\underline{\lambda}=\min\{\lambda_2,\,\lambda_{h}\}$ and $\bar\lambda=\max\{\lambda_{N},\,\lambda_H\}$. We summarize the effectiveness of this algorithm as follows and provide its proof in Appendix. 

\begin{lem}\label{thm:exact}
	Suppose Assumption \ref{ass:graph}  holds. Let $\kappa_1\geq 1$ and $\kappa_2\geq \frac{6N\bar \lambda^4\kappa_1^2}{\underline{\lambda}^4}\max\{1,\,\frac{1}{\underline{\lambda}}\}$. Then along the trajectory of   \eqref{alg:full-exact}, $x_i(t)$ exponentially converges to $y^*$ as $t$ goes to $\infty$.  
\end{lem}


\begin{rem}
In this algorithm we add two parameters  $\kappa_1$ and  $\kappa_2$ to compensate the difference between $L$ and $L^\top$ so as to handle directed communication graphs. Similar ideas have also been exploited in \cite{gharesifard2014distributed,kia2015distributed}. Particularly, when graph $\G$ is undirected, this solver can be simplified as 
\begin{align}\label{alg:full-exact-undirected}
	\begin{split}
	\dot{x}_i&=- {H}_i^\top [{H}_ix_i-z_i]- x_{oi}-v_{oi}\\
	\dot{v}_i&= x_{oi}
	\end{split}
\end{align}
which is consistent with existing solvers in \cite{liu2019arrow}. 
\end{rem}

However, in practice, we may only have noisy measurement $\tilde z_i$ instead of $z_i$. As  $\tilde z_i(t)$ can significantly deviate from its nominal value $z_i$, the performance of this algorithm cannot be guaranteed to produce the expected least squares solution. We will develop a novel identifier-based module to reject multiple sinusoidal disturbances even without knowing their frequencies, magnitudes, and initial phases. 

\subsection{Frequency identifier}

In this subsection, we put  $\tilde z_i$ into a parametric form and then develop local parameter identifiers for the spectrum of disturbance $\e_i$. This will facilitate us to solve the formulated problem using noisy measurements. 

Note that the disturbance $\e_i(t)$ under Assumption \ref{ass:disturbance} can be put into a generic form of $$\e_i(t)=\sum_{j=1}^{k_i} A_{ij} \sin(\omega_{ij} t+\phi_{ij})$$ for some constants $A_{i1},\,\dots,\,A_{ik_i}$ and $\phi_{i1},\,\dots,\,\phi_{ik_i}$. Thus we can rewrite the measurement equation \eqref{eq:main-noise}  as follows: 
\begin{align}\label{sys:full-measurement}
	\begin{split}
		\dot{\eta}_{i}&=S_{i}\eta_{i}\\
		\tilde z_i&=D_{i}\eta_{i}
	\end{split}
\end{align}
with $\eta_{i} \in \R^{2k_i+1}$ and   
\begin{align*}
S_{i}=&{\rm blkdiag}\left(0,\, \begin{bmatrix}0 & \omega_{i1}\\ -\omega_{i1}&0\end{bmatrix},\dots,\begin{bmatrix}0&\omega_{ik_i}\\
-\omega_{ik_i}&0\end{bmatrix}\right)    \\
D_{i}=&[\,1~ {1~0~1~\dots~1~0}\,]
\end{align*}
Here the initial condition ${\eta}_{i}(0)$ is determined by the constants $z_i$, $A_{i1},\,\dots,\,A_{ik_i}$, and $\phi_{i1},\,\dots,\,\phi_{ik_i}$. 
%
%

From here, we are tempted to develop an observer for  the internal state $\eta_i$ and separate the disturbance $\e_i(t)$ from the nominal value ${H}_i{y}$ in $\tilde z_i$ as that in \cite{tang2018iet,tang2018smc}. However, this procedure fails to be implemented without knowing the exact spectrum $\{\w_{i1},\,\dots,\,\w_{ik_i}\}$ of $\e_i(t)$.  Thus, we develop a parameter identifier to enable such an observer-based design to reject the considered disturbance as in \cite{xia2002global, sastry2011adaptive}.  

At first, one can practically verify that the pair $(D_i,\,S_i)$ is indeed observable by Lemma \ref{lem:observer} in Appendix. Then the observability matrix of system \eqref{sys:full-measurement} given by
\begin{align}
\Phi_i=\begin{bmatrix} D_i\\ D_iS_i\\ \vdots\\ D_iS_i^{2k_i}
\end{bmatrix} \in \R^{(2k_i+1)\times(2k_i+1)}
\end{align} 
is nonsingular.  By letting  $\eta_{i0}=\Phi_i \eta_i$, we have  
\begin{align}\label{sys:full-steady-state}
	\begin{split}
		\dot{\eta}_{i0}&= S_{i0}\eta_{i0}\\
		\tilde z_i&= D_{i0}\eta_{i0}
	\end{split}
\end{align}
where 
$$S_{i0}=\begin{bmatrix}
	{\bf 0}&I_{2k_i}\\
	0& {\bf p}_i
\end{bmatrix},\quad   D_{i0}=[1~~{\bf 0}_{2k_i}]$$ 
with vector ${\bf p}_i=-[\alpha_{i1}~0~\dots~\alpha_{ik_i}~0]\in \R^{2k_i}$ and constants $\alpha_{i1},\,\dots,\,\alpha_{ik_i}$ determined as the nonzero coefficients of the following polynomial:
$$p_i(s)\triangleq \prod\limits_{j=1}^{k_i}(s^2+w_{ij}^2)=s^{2k_i}+\alpha_{ik_i}s^{2k_i-2}+\dots+\alpha_{i2}s^{2}+\alpha_{i1}$$ 
From this we can also confirm that 
\begin{align*}
	\frac{{\rm d}^{(2k_i+1)}\,{\tilde z}_i}{{\rm d} \, t^{(2k_i+1)}}=-\alpha_{i1}\frac{{\rm d} \,\tilde z}{{\rm d }\,t}-\dots-\alpha_{ik_i}\frac{{\rm d }^{(2k_i-1)}\, {\tilde z}_i}{{\rm d }\, t^{(2k_i-1)}}
\end{align*}

Next, we choose two matrices $\tilde S_i\in \R^{(2k_1+1)\times (2k_i+1)}$ and $B_i\in \R^{(2k_i+1)\times 1}$:
$$\tilde S_{i}=\begin{bmatrix}
	{\bf 0}&I_{2k_i}\\
	-\tilde \beta_{i0}& \tilde {\bf p}_i
\end{bmatrix}, \quad  B_i=\begin{bmatrix}
{\bf 0}_{2k_i}\\1\end{bmatrix}$$
where $\tilde {\bf p}_i\triangleq -[\tilde \alpha_{i1}~\tilde \beta_{i1}~\dots~\tilde \alpha_{ik_i}~\tilde \beta_{ik_i}]$. The parameters $ \tilde\alpha_{ij}$ and $\tilde \beta_{ij}$ are selected such that the matrix $\tilde S_i$ is Hurwitz. 

Define a filtered output as follows:
\begin{align*}
	\hat z_i=[\tilde \beta_{i0}~ \tilde\alpha_{i1}-\alpha_{i1}~\tilde \beta_{i1}~\dots~\tilde \alpha_{ik_i}-\alpha_{ik_i}~\tilde \beta_{ik_i}]\tilde \eta_i
\end{align*} with $\tilde \eta_i\in \R^{2k_i+1}$ generated by 
\begin{align}\label{sys:filter}
	\begin{split}
		\dot{\tilde \eta}_{i}&=\tilde S_{i} \tilde \eta_{i}+ B_i \tilde z_i
	\end{split}
\end{align}
Following a similar procedure of Section 2.2 in \cite{sastry2011adaptive}, we confirm that  $\tilde z_i(t)=\hat z_i(t)+\bar{z}_i(t)$ holds with an exponentially vanishing error ${\bar z}_i(t)$ due to the initial conditions of $\tilde \eta_i(0)$.  

Let $\alpha_i=[\alpha_{i1}~\dots~\alpha_{ik_i}]^\top$ and $\hat \alpha_i=[\hat\alpha_{i1}~\dots~\hat\alpha_{ik_i}]^\top$ with $\hat \alpha_{ij}$ the estimate of $\alpha_{ij}$. The identifier output 
\begin{align*}
	\hat z_{o i}=[\tilde \beta_{i0}~\tilde \alpha_{i1}-\hat\alpha_{i1}~\tilde \beta_{i1}~\dots~\tilde \alpha_{ik_i}-\hat \alpha_{ik_i}~\beta_{ik_i}]\tilde \eta_i
\end{align*}
differs from the measurement output $\tilde z_i(t)$ by an identifier error 
\begin{align*}
	e_i=\hat z_{oi}-\tilde z_i	
\end{align*}
due to the initial conditions of $\tilde \eta_i(0)$ and the estimation error $\bar \alpha_i=\hat \alpha_i-\alpha_i$. To learn the vector ${\alpha}_i$, we introduce the following standard gradient algorithm:
\begin{align}\label{alg:identifier}
	\dot{\hat \alpha}_i={\ell} e_i \tilde \eta_{oi}	
\end{align}
with a learning rate $\ell>0$ and ${\tilde \eta}_{oi}=[{\tilde \eta}_{i2},\,\dots,\,{\tilde \eta}_{i(2k_i)}]^\top \in\R^{k_i}$. Having these estimates, we can get back to the polynomial $p_i(s)$ and  uniquely determine a time-varying and convergent estimate $\hat \w_{ij}(t)$ for $\w_{ij}$ analytically or numerically for the following analysis. 
 
\begin{rem}
	In the above identifier \eqref{alg:identifier}, we utilize the standard gradient algorithm to learn the unknown parameters. There are many alternative rules that can serve the same purpose, e.g., the normalized gradient algorithm of the form $\dot{\hat \alpha}_i={\ell} e_i \frac{\tilde \eta_{oi}}{1+\nu ||\tilde \eta_{oi}||^2}$ with a weight $\nu>0$. 
\end{rem}

\subsection{Distributed solvers with disturbance rejection}

With the parameter identifier \eqref{alg:identifier}, we substitute $\w_{ij}$ in matrix $S_i$ by its corresponding estimate $\hat \w_{ij}(t)$ and let 
\begin{align*}
\check S_{i}(t)=&{\rm diag}\left(0,\, \begin{bmatrix}0 & \hat  \omega_{i1}\\ -\hat \omega_{i1}&0\end{bmatrix},\dots,\begin{bmatrix}0&\hat \omega_{ik_i}\\
-\hat \w_{ik_i}&0\end{bmatrix}\right) 
\end{align*}
It can be verified that the two matrices $\check S_i$ and $\tilde S_i$ have no common eigenvalues under Assumption \ref{ass:graph}.  

Consider the following Sylvester equation:
	\begin{align}\label{eq:sylvester-est}
		\check T_i\check S_i-\tilde S_i \check T_i=B_iD_i
	\end{align}
The observability of $(D_i,\,\check S_i)$ and controllability of $(\tilde S_i,\, B_i)$ imply that  equation \eqref{eq:sylvester-est} must have a unique and invertible solution $\check T_i$ according to Theorem 2 in \cite{luenberger1965invertible}.  With this matrix $\check T_i$, we present the following distributed solver to determine the expected least squares solution irrespective of the considered measurement disturbance: 
\begin{align}\label{alg:full-adaptive}
		\dot{x}_i&=-\kappa_1 {H}_i^\top [{H}_ix_i-\tilde z_i +D _{i0} \check T_i^{-1}\tilde \eta_{i}]-(\kappa_1+\kappa_2) x_{oi}-v_{oi} \nonumber\\
		\dot{v}_i&=\kappa_1\kappa_2 x_{oi}\nonumber\\
		\dot{\tilde \eta}_{i}&=\tilde S_{i} \tilde \eta_{i}+ B_i \tilde z_i \\
		\dot{\hat \alpha}_i&=\ell  e_i \tilde \eta_{oi} \nonumber 
\end{align}
where $D_{i0}\triangleq [0~1~0~\dots~1~0]\in \R^{1\times (2k_i+1)}$ and the rest of these parameters are defined as above.  

This solver consists of two parts: a primal-dual dynamics to solve the nominal equations \eqref{eq:main} and an identifier-based compensator to estimate and cancel the  measurement disturbance. It is distributed in the sense that each node only uses its own and exchanged information from its neighbors in the network. 

\begin{rem}
	It is remarkable that different kinds of error-based adaptive rules have been developed in output regulation literature to reject such type of external disturbances with unknown or unknown frequencies, e.g., \cite{nikiforov2001nonlinear, xu2010robust}. Similar designs have been extended to solve some distributed optimization problem with external disturbances in \cite{wangxinghu2017, tang2018smc}. In comparison, the regulation error is not well-defined in our setting since the disturbance is acting on the measurement channel and the derived estimate may not converge to $y^*$. Thus we introduce a filter \eqref{sys:filter} and a filtered output for feedback. This converts the original disturbance rejection problem into an identification problem, and enables us 
	an identifier-based approach to reject the measurement disturbance. 
	Moreover, we can increase the order of the filter \eqref{sys:filter}  to improve the accuracy of our estimate $x_i$ relative to $y^*$ under noisy environments.   
\end{rem}

\section{Solvability analysis}\label{sec:sol}

In this section, we show the solvability of our problem via the preceding distributed solver \eqref{alg:full-adaptive}.  We first show the convergence of our solver towards $y^*$ irrespective of the measurement disturbance, and then verify its robustness with respect to possible unstructured but bounded disturbance. 

Here is the main result of this paper.
 
\begin{thm}\label{thm:unknown}
	Suppose Assumptions \ref{ass:graph} and \ref{ass:disturbance} hold. Then, for any initial point $x_i(0)$,\,$v_i(0)$,\, $\tilde\eta_i(0)$, and $\hat \alpha_i(0)$,  along the trajectory of \eqref{alg:full-adaptive},  ${x}_i(t)$ exponentially converges to the least squares solution ${y}^*$ to the system of linear equations \eqref{eq:main} as $t$ goes to $\infty$.  
\end{thm}
\pb To prove this theorem, we first put our solver into a compact form. Consider the Sylvester equation: $T_iS_i-\tilde S_i T_i=B_iD_i$. It has a unique and invertible solution $T_i$ by Theorem 2 in \cite{luenberger1965invertible}. Let $\bar \eta=[\bar \eta_1^\top,\,\dots,\,\bar \eta_N^\top]^\top$  and $\bar \alpha=[\bar \alpha_1^\top,\,\dots,\,\bar \alpha_N^\top]^\top$ with ${\bar \eta}_i=\tilde \eta_i-T_i\eta_i$. Then we obtain that:
\begin{align}\label{alg:full-adaptive-compact}
	\dot{x}&=-\kappa_1\check{H}^\top (\check{H} x - z)-(\kappa_1+\kappa_2)(L \otimes {I_m}){x}-(L \otimes {I_m})v  \nonumber\\
	&-\kappa_1\check{H}^\top {\check D}_0 \Delta  \nonumber\\
	\dot{v}&=\kappa_1\kappa_2(L \otimes {I_m})  {x}\\
	\dot{\bar \eta}&=\tilde S {\bar \eta}\nonumber\\
	\dot{\bar \alpha}&=-\ell {\tilde \eta}_{o}  {\tilde \eta}_o^\top \bar\alpha -\ell \bar {z} {\tilde \eta}_{o}    \nonumber
\end{align}
where  $\tilde \eta_{o} \triangleq [\tilde \eta_{o1}^\top,\,\dots,\, \tilde \eta_{oN}^\top]^\top$, $\bar {  z}=\mbox{diag}(\bar { z}_1,\,\dots,\,\bar { z}_N)$, ${\check D}_0=\mbox{blkdiag}(D_{1 0},\,\dots,\,D_{N0})$,  $\tilde S=\mbox{blkdiag}(\tilde S_{1},\,\dots,\,\tilde S_{N})$, and $\Delta=\eta-{\check T}^{-1}\tilde \eta$ with ${\check T}=\mbox{blkdiag}(\check T_{1},\,\dots,\,\check T_{N})$ and $\eta=[\eta_1^\top,\,\dots,\,\eta_N^\top ]^\top$.

To establish the convergence of $x$, we only have to consider the first two subsystems of \eqref{alg:full-adaptive-compact}. Following a similar procedure in the proof of Lemma \ref{thm:exact}, it is sufficient for us to focus on the stability issue of the following subsystem:
\begin{align}\label{sys:reduce-error-perturbed}
\dot{\bar{x}}&=-\kappa_1\check{H}^\top \check{H} \bar{x} -(\kappa_1+\kappa_2)(L \otimes {I_m}){\bar {x}}-(L {R} \otimes {I_m}){\tilde {  v}} \nonumber\\
&-\kappa_1 \check{H}^\top {\check D}_0 \Delta  \\
\dot{\tilde {v}}&=\kappa_1\kappa_2({R}^\top L \otimes {I_m}) {\bar {x}}\nonumber
\end{align}
Note that this subsystem is globally exponentially stable at the origin when $\Delta\equiv 0$ by the proof of Lemma \ref{thm:exact}. We thus view $\Delta$ as a perturbation and establish its stability as a perturbed system.  To complete the proof we are going to show that this perturbation is exponentially vanishing when $t$ tends to $\infty$. 

For this purpose we rewrite the perturbation into two parts as $\Delta=-T^{-1}\bar \eta-({\check T}^{-1}-T^{-1})\tilde \eta$ with $T=\mbox{blkdiag}(T_{1},\,\dots,\, T_{N})$. We shall show that each part  exponentially converges to $0$ when $t$ tends to $\infty$. 

On the one hand, by the choice of $\tilde S_i$, matrix $\tilde S$ is naturally Hurwitz. As a result, the $\bar \eta$-subsystem of \eqref{alg:full-adaptive-compact} is globally exponentially stable at the origin. That is, the trajectory of $\bar \eta$ will exponentially converge to $0$ when $t$ tends to $\infty$. This combined with the boundedness of $T^{-1}$ implies that the term $-T^{-1}\bar \eta$ exponentially vanishes when $t$ goes to $\infty$. 

On the other hand,  we focus on the second part. Let us start with the parameter identifier.  By Lemma 2.6.7 in \cite{sastry2011adaptive}, we obtain the persistence of excitation of each component of $\tilde \eta_{oi}$. Then the exponential convergence of $\bar \alpha$ can be ensured by Theorem 2.5.3 in \cite{sastry2011adaptive} when $\bar { z}(t)\equiv {\bf 0}$. Viewing $\bar { z}$ as an exponentially vanishing perturbation and using the boundedness of ${\tilde \eta}_{o}$, we can finally conclude the exponential convergence of $\bar \alpha$ towards $0$ by Corollary 9.1 in \cite{khalil2002nonlinear}.  Thus $||\check S_i(t)-S_i||$ will exponentially converge to $0$ as $t$ tends to $\infty$. Meanwhile, we can check that $(\check T_i-T_i)S_i-\tilde S_i(\check T_i-T_i)=-\check T_i(\check S_i-S_i)$. These two facts imply that $||\check T_i(t)-T_i||$ and $||\check T(t)-T||$ are exponentially vanishing as $t$ tends to $\infty$. We then recall the decomposition that $A^{-1}-B^{-1}=A^{-1}(B-A)B^{-1}$ for any two nonsingular matrices with compatible dimensions \cite{stewart1977perturbation}, and conclude that both $||\check T^{-1}-T^{-1}||$ and $({\check T}^{-1}-T^{-1})\tilde \eta $ will exponentially converge  to $0$ as $t$ tends to $\infty$. 

Overall,  $\Delta$ is indeed exponentially vanishing as the sum of two exponentially vanishing terms.  According to Corollary 9.1 in \cite{khalil2002nonlinear}, we can conclude the exponential convergence of \eqref{sys:reduce-error-perturbed} at the origin. As such, the trajectory of $x_i(t)$ will exponentially converge to $y^*$ as $t$ tends to $\infty$ under the proposed distributed solver \eqref{alg:full-adaptive}. The proof is thus complete.  
\pe


Next, we move on to the robustness issue of our solver with respect to unstructured disturbances. Consider the following measurement equation in this case: 
\begin{align*}
\tilde z_i={H}_i{y}+\e_i(t)+w_i(t)
\end{align*}
with $\e_i(t)$ the modeled disturbance fulfilling Assumption \ref{ass:disturbance} to be rejected and $w_i(t)$ some unstructured disturbance. Since $w(t)=[w_1(t),\,\dots,\,w_N(t)]^\top$ is often with high frequencies but quite small, we view it as a nonvanishing perturbation and aim to attenuate its effect on the convergence error $x_i(t)-y^*$ in the sense of finite-gain stability \cite{khalil2002nonlinear}.

\begin{thm}\label{thm:robust}
	Suppose Assumptions \ref{ass:graph} and \ref{ass:disturbance} hold.  Then, for any initial point $x_i(0)$,\,$v_i(0)$,\, $\tilde \eta_i(0)$, and $\hat \alpha_i(0)$,  along the trajectory of \eqref{alg:full-adaptive},  the inequality 
	\begin{align}\label{eq:finite-gain}
		\int_{0}^\infty |x_i(t)-{y}^*|^2 {\rm d}t \leq \gamma \int_{0}^{\infty}|w(t)|^2{\rm d }t+\delta 
	\end{align} 
	holds for some constants $\gamma>0$ and $\delta>0$.
\end{thm}
\pb  In this case, the composite system is still of the form  \eqref{alg:full-adaptive} and can be put into a compact form as follows:
\begin{align}\label{alg:full-adaptive-compact-robust}
\dot{x}&=-\kappa_1\check{H}^\top (\check{H} x - z-w) -(\kappa_1+\kappa_2)(L \otimes {I_m}){x}\nonumber \\
&-(L \otimes {I_m})v - \kappa_1 \check{H}^\top {\check D}_0 \Delta \nonumber\\
\dot{v}&=\kappa_1\kappa_2(L \otimes {I_m})  {x}  \\
\dot{\bar \eta}&=\tilde S {\bar \eta} \nonumber\\
\dot{\bar \alpha}&=-\ell {\tilde \eta}_{o}  {\tilde \eta}_o^\top \bar\alpha - \ell \bar { z} {\tilde \eta}_{o}    \nonumber
\end{align}
with an extra term related to $w$. We can take a similar coordinate transformation and derive an error system as follows:
\begin{align}\label{sys:reduce-error-perturbed-robust}
\dot{\bar{x}}&=-\kappa_1\check{H}^\top \check{H} \bar{x} -(\kappa_1+\kappa_2)(L \otimes {I_m}){\bar {x}}-(L {R} \otimes {I_m}){\tilde {v}}\nonumber\\
&-\kappa_1\check{H}^\top {\check D}_0 \Delta +\kappa_1\check{H}^\top w \nonumber \\
\dot{\tilde {v}}&=\kappa_1\kappa_2({R}^\top L \otimes {I_m}) {\bar {x}} \nonumber\\
\dot{\bar \eta}&=\tilde S {\bar \eta} \\
\dot{\bar \alpha}&=-\ell {\tilde \eta}_{o}  {\tilde \eta}_o^\top \bar\alpha + \ell \bar {\hat z} {\tilde \eta}_{o}  \nonumber 
\end{align}
Note that the error system consists of a cascade connection of the first two and last two subsystems. Recalling the global exponential stability of the $\bar\eta $-subsystem, $\bar \alpha $-subsystem, and the system \eqref{sys:reduce-error-perturbed} showed in Theorem \ref{thm:unknown}, we conclude that the above overall error system is also globally exponentially stable at the origin when $w=0$. Then we can directly recall Theorems 4.14 and 5.1 in \cite{khalil2002nonlinear} or follow the proofs to determine two constants to meet the requirement.  
\pe

In previous designs, we assume the multi-sinusoidal disturbance is totally unknown. When the spectrum of $\e_i(t)$ is known in prior, we can of course still use the solver \eqref{alg:full-adaptive} to tackle the problem. Meanwhile, it can be further simplified by removing the identification part as follows:
\begin{align}\label{alg:full-MN}
\dot{x}_i&=-\kappa_1 {H}_i^\top [{H}_ix_i-\tilde z_i +D _{i0}T_i^{-1} \tilde \eta_{i}]-(\kappa_1+\kappa_2) x_{oi}-v_{oi}\nonumber\\
\dot{v}_i&=\kappa_1\kappa_2 x_{oi}\\
\dot{\tilde \eta}_{i}&=\tilde S_{i} \tilde \eta_{i}+ B_i \tilde z_i	\nonumber
\end{align}
In this case the $\tilde \eta_i$-subsystem reduces to a Luenberger observer to estimate the full internal state of \eqref{sys:full-measurement}.  A particular choice for $\tilde S_i$ and $B_i$ is $\tilde S_i=S_i-K_{i}D_i$ and $B_i=K_i$ with $K_i\in \R^{1\times (2k_i+1)}$ a gain matrix such that $\tilde S_i$ is Hurwitz as shown in \cite{zhang2022ccc}. The existence of such gain matrices is guaranteed by the observability of $(D_i,\,S_i)$ for each $i$.

Here is the theorem and we omit its proof to save space. 

\begin{thm}\label{thm:known}
	Suppose Assumptions \ref{ass:graph} and \ref{ass:disturbance} hold.  Then, for any initial point $x_i(0)$, $v_i(0)$, and $\tilde\eta_i(0)$,  along the trajectory of \eqref{alg:full-MN},  ${x}_i(t)$ exponentially converges to the least squares solution ${y}^*$ to the system of linear equations \eqref{eq:main} as $t$ goes to $\infty$ when $w_i(t)\equiv 0$ and the inequality  \eqref{eq:finite-gain} holds for some constants $\gamma>0$ and $\delta>0$ when $w_i(t)\neq 0$.
\end{thm}

\begin{rem}\label{rem:disturbance}
	Compared with most existing distributed solvers \cite{mou2015distributed, liu2017exponential,shi2017network,liu2019arrow, wang2019distributed, zeng2018distributed}, Theorems \ref{thm:unknown}--\ref{thm:known} enable us to solve the problem under a general class of nonvanishing measurement disturbances. With such algorithms, the multi-sinusoidal disturbance can be asymptotically rejected even its frequency, amplitude, and initial phase are completely unknown while the unstructured disturbance is attenuated and has limited impact on the algorithm performance. This makes our distributed solvers more favorable in practical and noisy scenarios.
\end{rem}

\section{Numerical example}\label{sec:simu}

In this section, we provide a numerical example to illustrate the performance of  the proposed algorithms.

Consider the problem \eqref{eq:main} with
\begin{align*}
H&=\begin{bmatrix}
	 0.0479 &0.7514 &   0.5931   & 0.1329\\
	0.0176  &  0.0724  &  0.2320  &  0.5721
	\end{bmatrix}^{\top}\\
z&=\begin{bmatrix}
		10&20&30&40
	\end{bmatrix}^{\top}
\end{align*}
The nominal linear algebraic equations have a unique least squares solution $y^*\approx \left[22.33~~65.85\right]^{\top}$. 

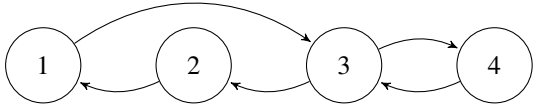
\begin{figure}
	\centering
	\begin{tikzpicture}[shorten >=1pt, node distance=2.cm, >=stealth',
		every state/.style ={circle, minimum width=1cm, minimum height=1cm}]
		\node[align=center,state](node1) {1};
		\node[align=center,state](node2)[right of=node1]{2};
		\node[align=center,state](node3)[right of=node2]{3};
		\node[align=center,state](node4)[right of=node3]{4};
		\path[->]  (node1) edge [bend left=35](node3)
		(node3) edge [bend left=25](node2)
		(node2) edge [bend left=25](node1)
		(node3) edge [bend left=25](node4)
		(node4) edge [bend left=25](node3)
		;
	\end{tikzpicture}
	\caption{Digraph $\mathcal{G}$ in our example. } \label{fig:graph}
\end{figure}

Suppose we have four nodes connected by a communication digraph as depicted in Fig.~\ref{fig:graph} with unity edge weights.  The measurement disturbance $\e_i(t)$ at node $i$ is a sinusoidal signal with frequency $\omega_i=0.5i$. We then resort to Theorems \ref{thm:unknown} and \ref{thm:known} to solve the problem. 

\begin{figure}
	\centering
	\includegraphics[width=0.40\textwidth]{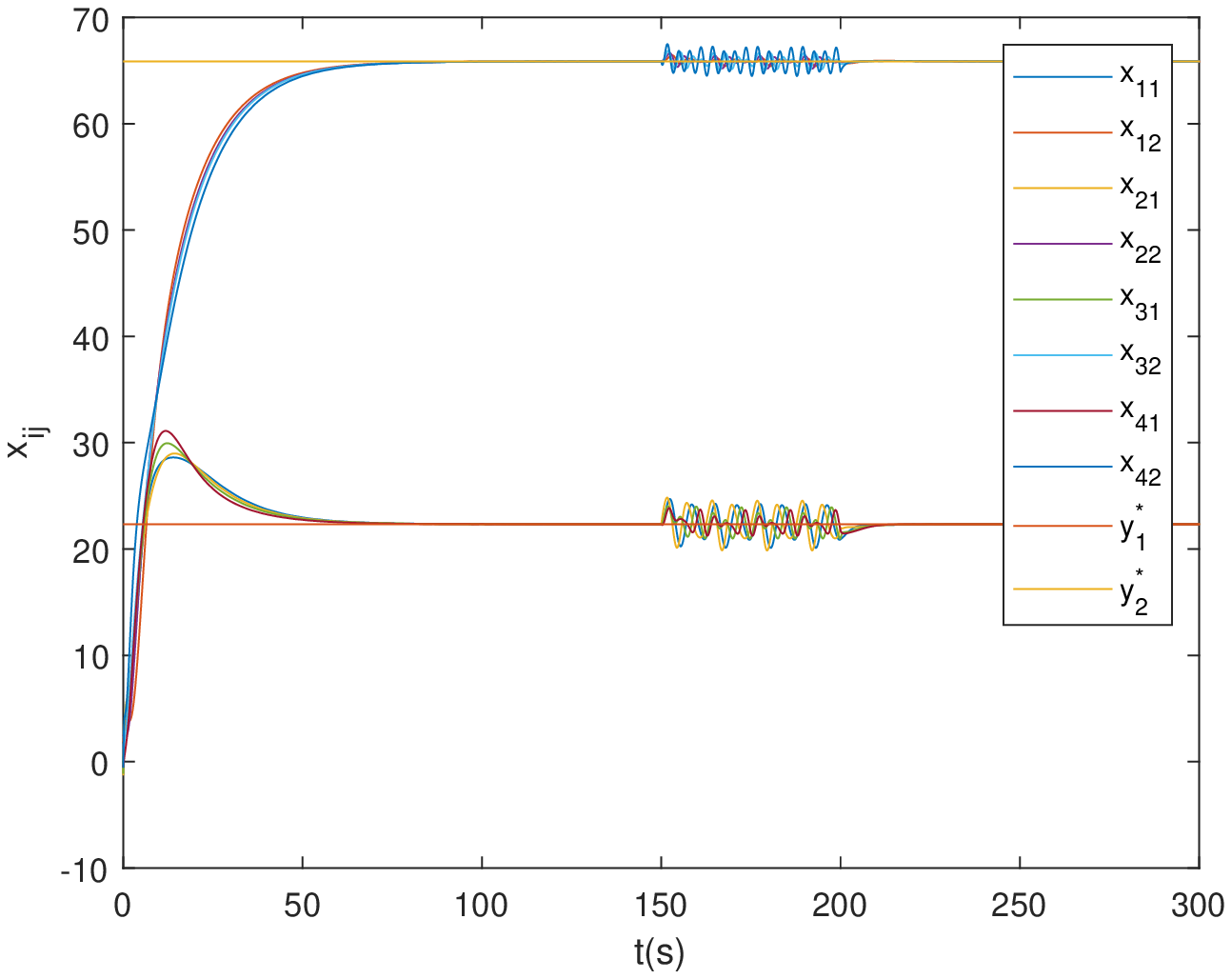}
	\caption{Profile of $x_i$ under algorithm \eqref{alg:full-MN}.} \label{fig:simu1}
\end{figure}

\begin{figure}
	\centering
	\includegraphics[width=0.42\textwidth]{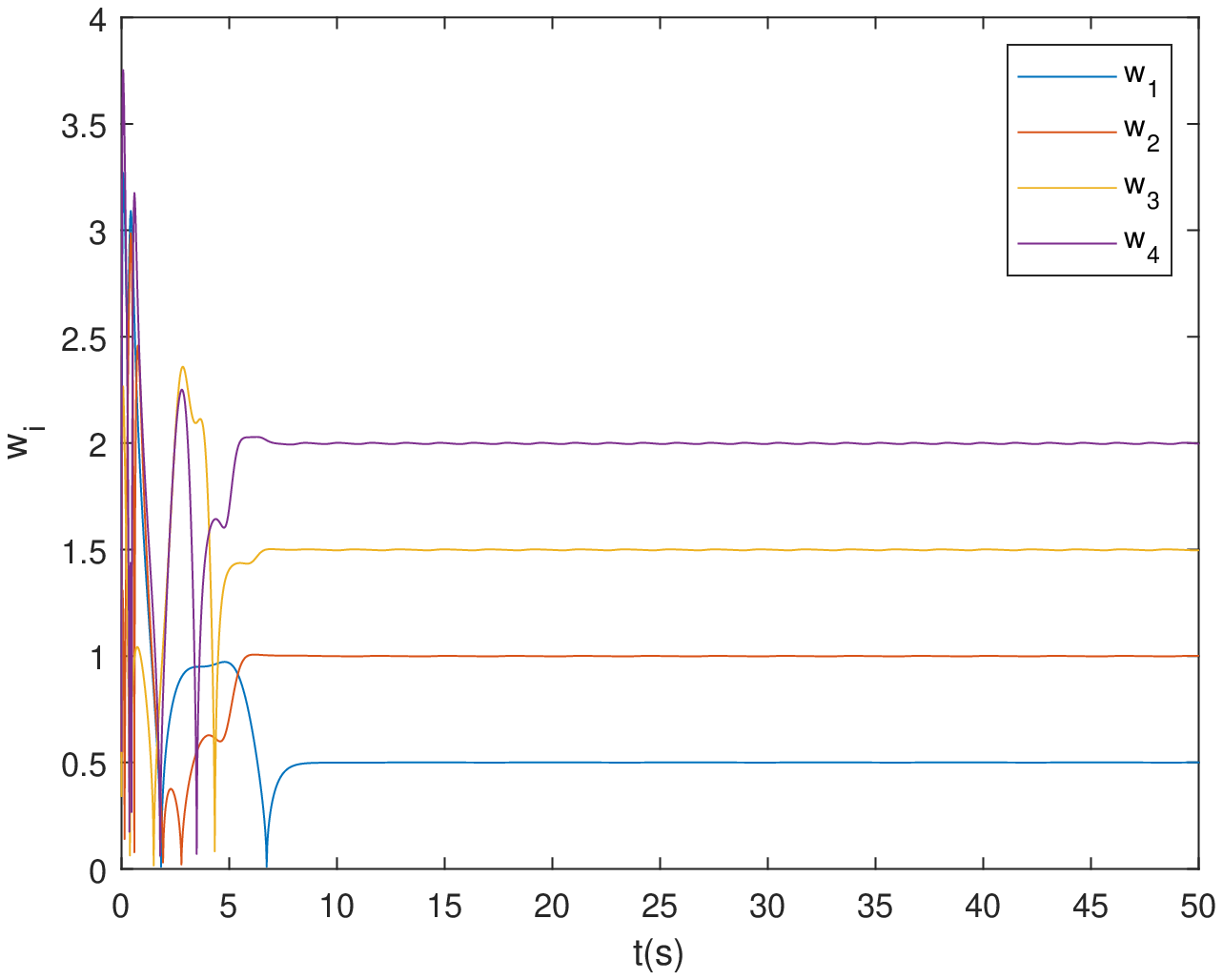}
	\caption{Estimate of frequencies under algorithm \eqref{alg:full-adaptive}.} \label{fig:simu2-w}
\end{figure}

In the simulation, we first consider the algorithm \eqref{alg:full-MN} with  $\tilde S_i=S_i-K_{i}D_i$ and $B_i=K_i$. Let us set
\begin{align*}
	K_1&=[24~-18~21.5]^\top,~~K_2=[6~0~10]^\top\\
	K_3&=[2.67~3.33~5.83]^\top, ~~K_4=[1.5~4.5~3.5]^\top
\end{align*}
All initial conditions are randomly chosen. The simulation result is shown in Fig.~\ref{fig:simu1} with $\kappa_1=\kappa_2=1$.  We shut down the disturbance rejection part at $t=150{\rm s}$ and reopen it after $t=200{\rm s}$.   The trajectory of $x_{i}$ is observed to converge to $y^*$ quickly at first and then deviate from the solution due to the measurement disturbance. After we reopen the disturbance rejection part, the convergence of $x_i(t)$ towards $y^*$ is quickly recovered as we expected.

\begin{figure}
	\centering
	\includegraphics[width=0.42\textwidth]{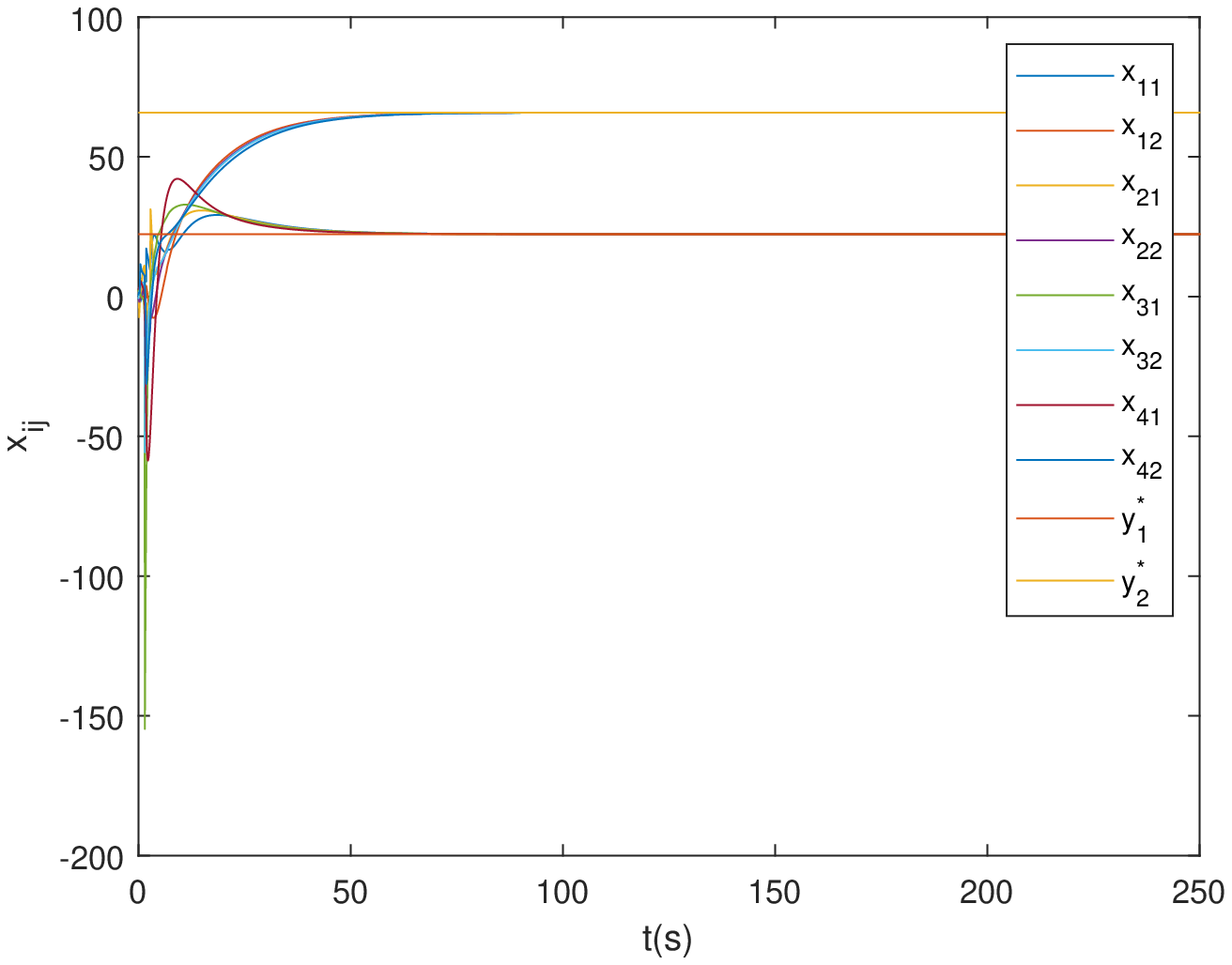}
	\caption{Estimate of  $x_i$ under algorithm \eqref{alg:full-adaptive}.} \label{fig:simu2-x}
\end{figure}

Next, we choose $\beta_{i0}=8$, $\alpha_{i1}=12$, $\beta_{i1}=6$ and use algorithm \eqref{alg:full-adaptive} to solve the problem. We employ the normalized gradient rule and set $\ell=30$ and $\nu=1$ for adaptation. The profile of estimated frequencies is reported in Fig.~\ref{fig:simu2-w}. We can find that the estimators indeed converge to the true frequencies. Then we list the trajectories of node states in Fig.~\ref{fig:simu2-x}. It is confirmed that the adaptive algorithm \eqref{alg:full-adaptive} works well in solving our problem even the frequencies of the measurement disturbances are totally unknown.

\begin{figure}
	\centering
	\includegraphics[width=0.42\textwidth]{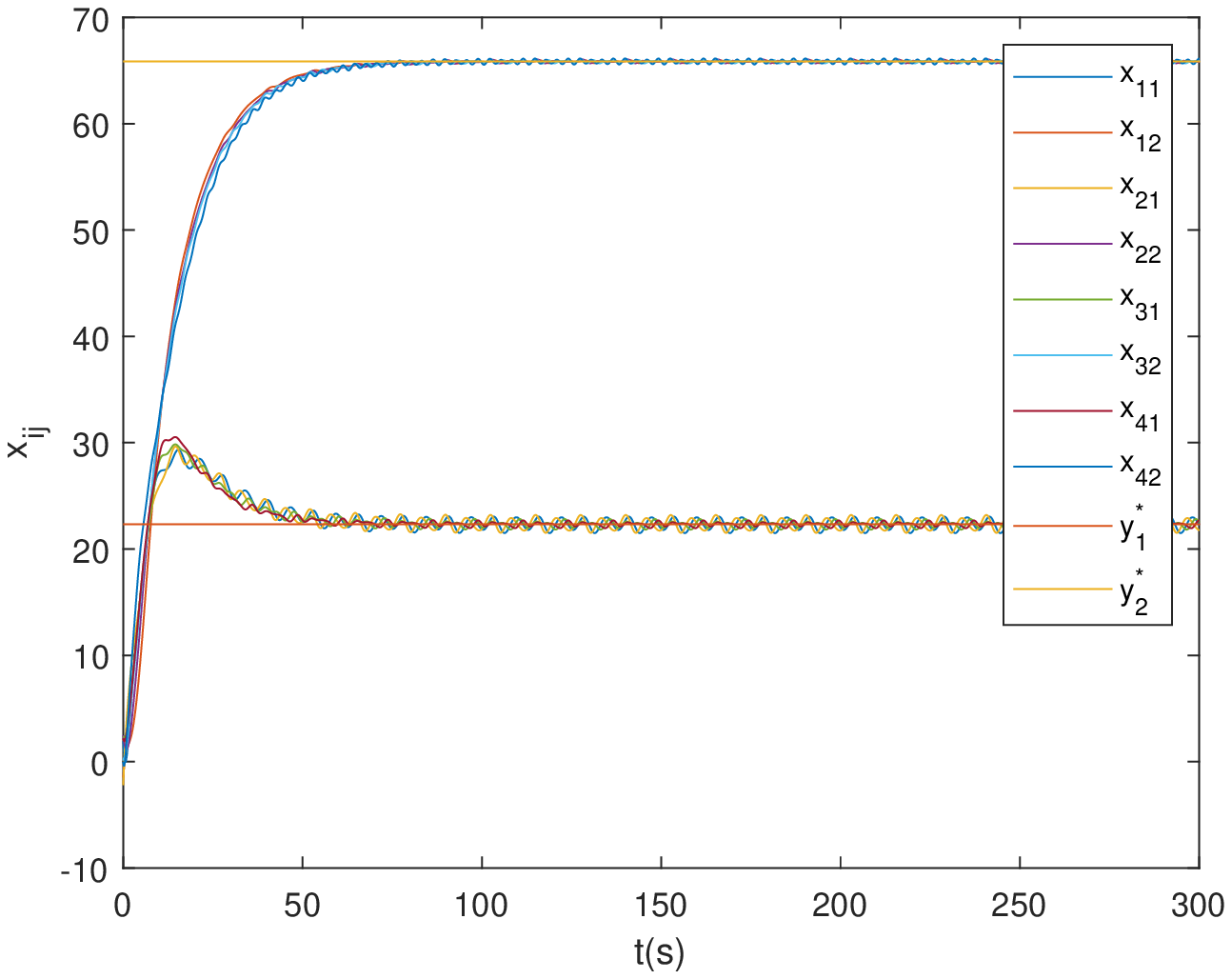}
	\caption{Estimate of   $x_i$ with washout filter $\frac{s}{s+d}$ and $d=0.4$.} \label{fig:simu2-wash}
\end{figure}

To make a comparison, we also use a washout filter to handle the measurement disturbance as in \cite{levine2011control}. The simulation result is given in Fig.~\ref{fig:simu2-wash}. The fluctuation of the trajectories is indeed suppressed in this case but cannot be totally rejected with the chosen parameter. These observations validate the efficacy of our preceding design to distributedly and exactly solve the least squares problem for given linear equations with disturbance rejection.

\section{Conclusion}\label{sec:con}
We have studied the distributed continuous-time algorithm for solving a system of linear algebraic equations under measurement disturbances. By inserting an (identifier-enhanced) observer to estimate the  disturbances, we have developed effective distributed algorithms for two different cases depending upon whether the frequencies of the multi-sinusoidal disturbances are known a prior or not. We have also shown the robustness of our designs against unstructured and bounded disturbances. In the future, we may investigate this problem for more general communication digraphs or extend such design to a fully distributed case.

\appendix  \label{app}

\subsection{Graph notion}

A weighted directed graph (digraph) is described by $\mathcal {G}=(\mathcal {N}, \mathcal {E}, \mathcal{A})$ with node set $\mathcal{N}=\{1,{\dots},N\}$ and edge set $\mathcal {E}$. $(i,\,j)\in \mathcal{E}$ denotes an edge from node $i$ to node $j$. The weighted adjacency matrix $\mathcal{A}=[a_{ij}]\in \mathbb{R}^{N\times N}$ is defined by $a_{ii}=0$ and $a_{ij}\geq 0$. Here $a_{ij}>0$ iff there is an edge $(j,\,i)$ in the digraph. Node $i$'s neighbor set is defined as $\mathcal{N}_i=\{j\mid (j,\, i)\in \mathcal{E} \}$. 
If there is a directed path between any two nodes, then the digraph is said to be strongly connected.  The in-degree and out-degree of node $i$ are defined by $d^{\mbox{in}}_i=\sum\nolimits_{j=1}^N a_{ij}$ and $d^{\mbox{out}}_i=\sum\nolimits_{j=1}^N a_{ji}$. A digraph is weight-balanced if $d^{\mbox{in}}_i=d^{\mbox{out}}_i$ for any $i\in \mathcal{N}$. The Laplacian of $\mathcal{G}$ is defined as $L\triangleq D^{\mbox{in}}-\mathcal{A}$ with $D^{\mbox{in}}=\mbox{diag}(d^{\mbox{in}}_1,\,\dots,\,d^{\mbox{in}}_N)$. 

\subsection{Proof of Lemma \ref{thm:exact}}

We first put the algorithm into a compact form:
\begin{align}\label{sys:compact-disturbance-free}
		\dot{x}&=-\kappa_1\check{H}^\top (\check{H} x - z)-(\kappa_1+\kappa_2)(L \otimes {I_m}){x}-(L \otimes {I_m})v \nonumber \\
		\dot{v}&=\kappa_1 \kappa_2(L \otimes {I_m}) {x}
\end{align}
with ${x}=[x_1~\cdots~x_N]^\top$, ${v}=[v_1~\cdots~v_N]^\top$, and $\check{H}=\mbox{blkdiag}({H}_1,\,\dots,\,{H}_N)$. 

To show the relationship between $y^*$ and the equilibrium point of \eqref{sys:compact-disturbance-free}, we set the righthand side of the above compact form to zero, and have some constant ${x}^*$ and ${v}^*$  such that 
\begin{align*}
	\kappa_1 \check{H}^\top (\check{H} {x}^*-z)+(L \otimes {I_m}){ {v}^*}={\bf 0}\\
	(L \otimes {I_m}) { {x}^*}={\bf 0}
\end{align*}
From the second equation, there should be some $s\in \R^m$ such ${x}^*={\bf 1}_N \otimes s$.  Using the fact that ${\bf 1}^\top L={\bf 0}$ under Assumption \ref{ass:graph}, we premultiply the first equation by ${\bf 1}^\top \otimes I_m$ and obtain  $\check{H}^\top (\check{H} {x}^*-z)={\bf 0}$, or equivalently, $H^\top (Hs-z)={\bf 0}$. Since the rank of $H$ is full, we have $s=y^*$ and ${x}^*={\bf 1}_N \otimes {y}^*$. Then, it suffices for us to prove the convergence of $x(t)$ towards $x^*$ along the trajectory of system \eqref{sys:compact-disturbance-free} when $t$ goes to $\infty$.  

For this purpose we let $\bar {x}={x}-{x}^*$, ${\tilde{ v}}_{0}=(\frac{{\bm 1}_N^\top}{\sqrt{N}} \otimes {I_m}){\bar{v}}$, and ${\tilde{ v}}=({R}^\top \otimes {I_m}){\bar{v}}$ with $\bar {v}={v}-{v}^*$. It follows that ${\tilde{ v}}_{0}\equiv 0$ and 
\begin{align}\label{sys:full-error-disturbance-free}
	\dot{\bar{x}}&=-\kappa_1\Xi -\kappa_2(L \otimes {I_m}){\bar {x}}-(LR \otimes {I_m}){\tilde {v}} \nonumber\\
	\dot{\tilde {v}}&=\kappa_1 \kappa_2(R^\top L \otimes {I_m}) {\bar {x}} 
\end{align}
with  $\Xi \triangleq [\check{H}^\top \check{H}+(L \otimes {I_m})]\bar x$. Here we claim that the following  inequality holds for any $\bar x$:
\begin{align}\label{eq:strong-monotone-disturbance-free}
	\bar x^\top \Xi \geq \frac{\underline{\lambda}}{N}\bar x^\top \bar x
\end{align} 
We consider two different cases. 
On the one hand, suppose $(L\otimes {I_m})\bar x\neq 0$. Under Assumption \ref{ass:graph}, we 
have 
\begin{align*}
	\bar x^\top \Xi\geq \bar x^\top (L \otimes {I_m}){\bar {x}}= \bar x^\top [\mbox{Sym}(L) \otimes {I_m}]{\bar {x}}\geq {\lambda}_{2}\bar x^\top \bar x
\end{align*}
On the other hand, when $(L\otimes {I_m})\bar x=0$, from the definition of $\bar x$, there exists some non-zero constant vector $s'\in \R^m$ such that $\bar x={\bf 1}_N \otimes s'$. In this case, it follows then
\begin{align*}
	\bar x^\top \Xi&=\bar x^\top \check{H}^\top \check{H} \bar{x}= {s'}^\top H^\top H s' \geq \frac{\lambda_{h}}{N} ||\bar x||^2
\end{align*}
Overall, the inequality \eqref{eq:strong-monotone-disturbance-free} is confirmed.

Next, we get back to system \eqref{sys:full-error-disturbance-free}.  To establish its stability we perform a coordinate transformation: $\hat x_u=(\frac{{\bf 1}^\top_{N}}{\sqrt{N}}\otimes I_m)\bar x$, $\hat x_l=(R^\top \otimes I_m)\bar x$, and $\hat v=\tilde v+ \kappa_1(R^\top \otimes I_m)\bar x$. This gives
\begin{align}\label{sys:reduce-error-disturbance-free-skew}
	\dot{\hat {x}}_u&=-\kappa_1(\frac{{\bf 1}^\top_{N}}{\sqrt{N}}\otimes I_m)\Xi   \nonumber\\
	\dot{\hat {x}}_l&=-\kappa_1(R^\top \otimes I_m)\Xi - (\kappa_2-\kappa_1)(R^\top L R\otimes I_m) \hat x_l\\ 
	&-(R^\top L {R} \otimes {I_m}){\hat {  v}} \nonumber\\
	\dot{\hat {v}}&=-\kappa_1 (R^\top L R \otimes I_m)\hat v+\kappa_1^2(R^\top L R\otimes I_m) \hat x_l \nonumber\\
	&-\kappa_1^2 (R^\top \otimes I_m)\Xi \nonumber
\end{align} 

Choose a positive definite Lyapunov function candidate as $V_{o}(\hat x_u,\,\hat x_l,\,\hat v)= \kappa_3{\hat x}^\top_u \hat x_u+ \kappa_3 \hat x_l^\top \hat x_l+\hat v^\top \hat v$ with $\kappa_3>1$ to be specified later. The derivative of $V_{o}$ along the trajectory of the system \eqref{sys:reduce-error-disturbance-free-skew} satisfies
\begin{align*}
	\dot{V}_{o}&=-2\kappa_1 \kappa_3 \bar x^\top  \Xi  -2 (\kappa_2-\kappa_1)\kappa_3\hat x_l(R^\top L R\otimes I_m) \hat x_l\\
	&- 2  \kappa_3 {\hat x}_l^\top (R^\top LR\otimes I_m)\hat v-2\kappa_1 {\hat v}^\top (R^\top L R\otimes I_m) \hat v\\
	&+2\kappa_1^2{\hat{v}}^\top (R^\top LR \otimes I_m)\hat x_l- 2\kappa_1^2{\hat{v}}^\top (R^\top \otimes I_m)\Xi\\
	&\leq -\frac{2\kappa_1 \kappa_3 \underline{\lambda}}{N} ||\bar x||^2-2\kappa_1 \underline{\lambda}||\hat v||^2+6\kappa_1^2\bar \lambda  ||\bar x||||\hat v||\\
	&-2 (\kappa_2-\kappa_1)\kappa_3 \underline{\lambda}||\hat x_l||^2+2 \kappa_3 \bar \lambda ||\hat x_l||\hat v||
\end{align*}
where we have used the fact that  $||\Xi||\leq 2\bar{\lambda}||\bar x||$ for any $\bar x$ and inequality   \eqref{eq:strong-monotone-disturbance-free}.  Completing the squares gives
\begin{align*}
	\dot{V}_{o}	&\leq -[\frac{2\kappa_1 \kappa_3 \underline{\lambda}}{N}-\frac{18\kappa_1^3\bar \lambda^2}{\underline{\lambda}}] ||\bar x||^2-\kappa_1 \underline{\lambda}||\hat v||^2\\
	&-[2 (\kappa_2-\kappa_1)\kappa_3 \underline{\lambda}-\frac{\kappa_3^2\bar \lambda^2}{\kappa_1\underline{\lambda}}]||\hat x_l||^2	
\end{align*}
Choose $\kappa_3=\frac{10N \bar \lambda^{2}\kappa_1^3}{ \underline{\lambda}^2}\max\{1,\,\frac{1}{\underline{\lambda}}\}$. Under the lemma conditions, one can obtain that 
\begin{align*}
	\dot{V}_{o}	&\leq -2\kappa_1^3||\bar x||^2-\kappa_1 \underline{\lambda}||\hat v||^2
\end{align*}
which implies the global exponential stability of system \eqref{sys:reduce-error-disturbance-free-skew} and \eqref{sys:full-error-disturbance-free} at the origin. The proof is thus complete. 

\subsection{A lemma on system observability}
Here is a lemma on the parallel connection's observability for two given observable systems. It assists us in checking the observability of system \eqref{sys:full-measurement}. 
\begin{lem}\label{lem:observer}
	Consider matrices $\Psi_1\in \R^{1\times n_1}$, $\Psi_2\in \R^{1\times n_2}$, $A_1\in \R^{n_1\times n_1}$, and $A_2\in \R^{n_2\times n_2}$ with $n_1,\,n_2\geq 1$. Assume that the pairs $(\Psi_1,\,A_1)$ and $(\Psi_2,\,A_2)$ are observable. Then the following pair 
	\begin{align*}
		\left( \begin{bmatrix}
			\Psi_1&\Psi_2
		\end{bmatrix},\quad  \begin{bmatrix}
			A_1&{\bf 0}\\
			{\bf 0}&A_2
		\end{bmatrix}\right)
	\end{align*}
	is  observable iff $A_1$ and $A_2$ have no common eigenvalues. 
\end{lem}
\pb The sufficiency is a direct consequence of the PBH test and we just prove the necessity by seeking a contradiction. 

We suppose $\lambda$ is a common eigenvalue of the matrices $A_1$ and $A_2$. By definition, there exist two vectors $\eta_1\in \R^{n_1}$ and $\eta_2\in \R^{n_2}$ such that  $A_i\eta_i=\lambda \eta_i$ hold for $i=1,\,2$. Due to the observability of $(\Psi_i,\,A_i)$, we must have $\Psi_i\eta_i\neq 0$. If not, $\lambda$ will correspond to an unobservable mode of the pair $(\Phi_i,\,A_i)$ by the PBH test. Without loss of generality, we assume  $\Psi_1\eta_1+\Psi_2\eta_2=0$. Otherwise, we can replace $\eta_1$ by $-\frac{\Psi_2\eta_2\eta_1}{\Psi_1\eta_1}$. Let $\eta=\mbox{col}(\eta_1,\,\eta_2)$. One can verify that 
\begin{align*}
	\begin{bmatrix}
		A_1&{\bf 0}\\
		{\bf 0}&A_2
	\end{bmatrix}\eta=\lambda\eta,\quad \begin{bmatrix}
		\Psi_1&\Psi_2
	\end{bmatrix}\eta =0
\end{align*} 
Thus we arrive at a contradiction with the observability premise by the PBH test. 
The proof is thus complete. 
\pe

\bibliographystyle{IEEEtran}
\bibliography{distr_eq}

\end{document}